# COMMUTING POLYNOMIALS AND POLYNOMIALS WITH SAME JULIA SET

Pau Atela
Jun Hu

April, 1995

ABSTRACT. It has been known since Julia that polynomials commuting under composition have the same Julia set. More recently in the works of Baker and Eremenko, Fernández, and Beardon, results were given on the converse question: When do two polynomials have the same Julia set? We give a complete answer to this question and show the exact relation between the two problems of polynomials with the same Julia set and commuting pairs.

§1 *Introduction and statement of results.*

Let $f, g \colon \hat{\mathbb{C}} \to \hat{\mathbb{C}}$ be two non-linear complex polynomial maps on the Riemann sphere with Julia sets $J_f, J_g$ and filled-in Julia sets $K_f, K_g$. Recall that $K_f$ is the set of bounded orbits of $f$ and $J_f = \partial K_f$. Julia [J] showed that if two rational functions commute ($f \circ g = g \circ f$), then necessarily $J_f = J_g$. For polynomials this is in fact easy to see:

If $\{f^i(z)\}$ is bounded, then $\{g(f^i(z)) = f^i(g(z))\}$ is also bounded. Therefore $g(K_f) \subset K_f$, and consequently $g^i(K_f) \subset K_f \ \ \forall i \geq 0$. This shows that $K_f \subset K_g$. Similarly $K_g \subset K_f$, so $J_f = J_g$, and the proof is complete.

Several authors (see [B-E], [F], [Bea1], [Bea2]) have considered the converse question: when do two polynomials have the same Julia set? We will summarize these previous results after introducing some notation. We have structured the presentation so as to emphasize the parallelism between the two problems of polynomials with same Julia set and commuting pairs.

We will denote by $S^1 \subset \mathbb{C}$ the unit circle. The composition $f \circ g$ will be written as $fg$, and the $i$–fold composition $f \circ \cdots \circ f$ as $f^i$. For a fixed polynomial $f$ of degree $n \geq 2$, we define

$\mathcal{S} = \{\, g \mid J_g = J_f \,\}$, the set of polynomials with same Julia set as $f$;

$\mathcal{C} = \{\, g \mid fg = gf \,\}$, the set of polynomials that commute with $f$; and

$\Sigma$ as the group of rotational symmetries of $J_f$.

We also define (for $m \geq 2$)

$\mathcal{S}_m = \{\, g \mid \deg(g) = m \text{ and } J_g = J_f \,\}$, and

$\mathcal{C}_m = \{\, g \mid \deg(g) = m \text{ and } fg = gf \,\}$.

Notice that $\mathcal{S} = \cup \mathcal{S}_m$, and $\mathcal{C} = \cup \mathcal{C}_m$.

Both author's research at MSRI was partially supported by NSF grant #DMS–9022140. Atela's research was also partially supported by NSF grant #DMS–9316919.

**Definition.** A polynomial $\sum c_i z^i$ of degree $d$ is called *centered* if $c_{d-1} = 0$.

It is easy to see that a polynomial $\tilde{f}$ is always conjugate to a centered one $f$ by a simple translation. That is, $\tilde{f} = L^{-1} f L$ with $L(z) = z + a$, for some constant $a \in \mathbb{C}$. The Julia sets are just translations of each other: $L(J_{\tilde{f}}) = J_f$.

For a centered polynomial $f$ the rotational symmetries are rotations around the origin. We can therefore write $\Sigma = \{ c \in \mathbb{C} \mid cJ = J \}$. Notice that necessarily $|c| = 1$. Clearly $\Sigma$ is a subgroup of $S^1$ and, since $J$ is closed, either $\Sigma = S^1$ or $\Sigma = \{l^{th} \text{ roots of unity}\}$ for some $l$. For a study on rotational symmetries of Julia sets of a polynomial we refer the reader to [Bea1].

*Remark.* Julia sets might have other than rotational symmetries. For instance, symmetry with respect to a line.

With this notation, previous results are:

- Julia [J, 1922]:  $\mathcal{C} \subset \mathcal{S}$.
- Baker and Eremenko [B-E, 1987]:  $\mathcal{C} = \mathcal{S}$ when $\Sigma$ is trivial.
- Fernández [F, 1989]:  If $f, g$ same degree and same leading coefficient, $J_f = J_g \implies f = g$.
- Beardon [Bea1, 1990]:  $\mathcal{S} = \{ g \mid gf = \sigma f g \text{ for some } \sigma \in \Sigma \}$.
- Beardon [Bea2, 1992]:  $\{ g \mid \deg(f) = \deg(g), J_f = J_g \} = \{ \sigma f \mid \sigma \in \Sigma \}$.

These papers contain no examples; our initial task, then, was to find concrete examples. Consider the following simple

**Examples.**
  (1) Consider the family $\{ z^n \mid n \geq 2 \}$. Any two members indeed commute and $J = S^1$.
  (2) Consider the family of Tchebycheff polynomials $\{ T_n(z) \mid n \geq 2 \}$, where $T_n(z) = \cos(n \arccos z)$. One easily checks that $T_n(T_m(z)) = T_{nm}(z)$, and so any two commute. The common Julia set is $J = [-1, 1]$.
  (3) Fix a polynomial $f$ and consider the family of iterates $\{ f^n \}$. Obviously they all commute and have the same Julia set.

Our first result shows that in fact these are, essentially, the only examples:

**Theorem 1.** *Let $f$ and $g$ be two centered polynomials of degrees $n, m \geq 2$ with same Julia set $J$. Then we have one of the following three cases:*
  (1) *If $J$ is a circle, then $\exists L$ linear map and $\sigma \in \Sigma = S^1$ such that $f = L z^n L^{-1}$ and $g = L \sigma z^m L^{-1}$.*
  (2) *If $J$ is an interval, then $\exists L$ linear map and $\sigma \in \Sigma = \{-1, 1\}$ such that $f = L T_n L^{-1}$, $g = L \sigma T_m L^{-1}$, where $T_n, T_m$ are Tchebycheff polynomials.*
  (3) *In all other cases, there is a polynomial $R$ and positive integers $q, Q$ such that*
  $$f = \epsilon_1 R^q, \qquad g = \epsilon_2 R^Q, \qquad \text{with } \epsilon_1, \epsilon_2 \in \Sigma.$$

The two problems—polynomials with same Julia set and commuting polynomials—are closely related. The next theorems give a complete answer showing exactly how.

For polynomials related to the examples (1) and (2), we have that $\mathcal{S}_m \neq \emptyset \neq \mathcal{C}_m \ \forall m \geq 2$. However, in general we have $\mathcal{S}_m = \mathcal{C}_m = \emptyset$ for most $m$.

**Definition.** In this paper, a polynomial $f$ will be called *minimal* if it is centered and $f \neq \sigma R^i$ with $\deg(R) < \deg(f)$, $\sigma \in \Sigma$ and $i > 1$.

**Theorem 2.** *If $f$ is a minimal polynomial of degree $n \geq 2$ and its Julia set is neither a circle nor an interval, then:*

(1) $\mathcal{S}_m \neq \emptyset \iff m = n^i$ *for some* $i \geq 1$;
(2) $\mathcal{S}_{n^i} = \Sigma f^i$; *and*
(3) $\mathcal{S} = \Sigma \mathcal{C}$.

From theorems 1 and 2, we immediately have

**Theorem 3.** *Let $f$ be a polynomial of degree $n$ with Julia set $J$.*

(1) *If $J$ is a circle or an interval, $\mathcal{S} = \Sigma \mathcal{C}$.*
(2) *Otherwise, if $f$ is minimal, $\mathcal{S} = \Sigma \mathcal{C}$.*

*Remark.* If $f$ is not minimal we don't necessarily have $\mathcal{S} = \Sigma \mathcal{C}$. For instance, let $R = z^2(z^3 + 1)$ and take $f = e^{\frac{2\pi i}{3}} R^2$, $g = R$; as we will see in the next section, we have $\Sigma = \{\, \sigma \mid \sigma^3 = 1 \,\} = \{1, e^{\frac{2\pi i}{3}}, e^{\frac{4\pi i}{3}}\}$, $g \in \mathcal{S}$ but $g \notin \Sigma \mathcal{C}$. That is, $J_g = J_f$ but $\sigma g$ does not commute with $f$ for any $\sigma \in \Sigma$.

*Remark.* The question for rational maps and meromorphic functions is more difficult. R. Devaney pointed out to us this simplest example:

The map $z \mapsto \frac{z^2}{1+2z}$ (which is $z \mapsto z^2$ conjugated by $w \mapsto \frac{1}{w-1}$) has $J = \mathbb{R}$ and so does $z \mapsto \tan z$.

§2 *The Boettcher function.*

It is well known that for a polynomial $f = az^n + \ldots$ of degree $n$ the point at infinity is super-attracting and so there is a function $\phi$ (called the Boettcher function), analytic in a neighbourhood of infinity, that satisfies

$$\phi f \phi^{-1} = az^n,$$

with power series expansion

$$\phi(z) = z + b_0 + \frac{b_1}{z} + \frac{b_2}{z^2} + \ldots$$

(often $\phi$ is required to satisfy $\phi f \phi^{-1} = z^n$.)

When the Julia set $J$ is connected, the extension of $\phi$ turns out to be the Riemann map of the complement of $K_f$. So it is natural to expect that if $g$ is another polynomial with $J_g = J_f$, then $g$ and $f$ will have the same Boettcher function. For details on the following see, for example, [B-E].

**Proposition 1.** *If $f(z) = az^n + \ldots$ and $g(z) = bz^m + \ldots$ are two polynomials of degrees $n, m \geq 2$ with same Julia set $J$, there is a unique analytic function $\phi$ in a neighbourhood of infinity such that*

$$\phi f = L_f \phi, \qquad \phi g = L_g \phi,$$

with power series expansion $\phi(z) = z + b_0 + \frac{b_1}{z} + \frac{b_2}{z^2} + \ldots$, where $L_f(z) = az^n$ and $L_g(z) = bz^m$.

*Sketch of proof.*
Let $\phi, \psi$ be the Boettcher functions of $f$ and $g$, respectively. If the Julia set $J$ were connected, no critical points escape to infinity and $\phi$ can be extended to an analytic homeomorphism $\hat{\mathbb{C}} \setminus K_f \to \hat{\mathbb{C}} \setminus D_1$, where $D_1$ denotes the closed unit disc. The equality would then follow from Schwarz's Lemma. However, in general one considers the function $G(z) = \log|\phi| + (\frac{1}{n-1})\log|a|$ (and the corresponding one for the map g). It is well known that this is the Green function of $\hat{\mathbb{C}} \setminus K_f$ with pole at $\infty$. As $K_f = K_g$, by unicity of the Green function we have

$$\log|\phi(z)| + \frac{1}{n-1}\log|a| = \log|\psi(z)| + \frac{1}{m-1}\log|b|.$$

One deduces that $|\phi(z)| = |\psi(z)|$, and since $\phi(z) \to z, \psi(z) \to z$, as $z \to \infty$, we must have $\phi(z) \equiv \psi(z)$. □

By equating coefficients of $z^{n-1}$ one gets

**Lemma 2.** *$f$ is centered iff $b_0 = 0$. Therefore, if $J_f = J_g$, $f$ is centered iff $g$ is centered.*

**Lemma 3.** *$J = J_f = J_g \implies J_{fg} = J_{gf} = J$.*

*Proof.* Let $K, K_{fg}$ denote the filled-in Julia sets of $f$ (and $g$) and $fg$ respectively, and $F, F_{fg}$ the Fatou sets. As $f(K) = K$ and $g(K) = K$, then $(f \circ g)^{\circ i}(K) = K \; \forall i \geq 0$. Therefore $K \subset K_{fg}$. On the other hand we also have $(f \circ g)^{\circ i}(F) = F \; \forall i \geq 0$, so by Montel's theorem the family $\{(f \circ g)^{\circ i}\}$ is normal in all of $F$, so $F \subset F_{fg}$. This implies that $F^\infty \subset F_{fg}^\infty$ (these being the connected components of the Fatou sets at $\infty$), and therefore we have for the complement sets $K \supset K_{fg}$. □

*Observation.* If $P$ is a polynomial of degree $p$, then expanding around infinity gives

$$\tfrac{1}{P(z)} = O(\tfrac{1}{z})^p.$$

**Lemma 4.** *Let $c \in \mathbb{C}$, $\phi$ an analytic function in a neighbourhood of infinity with $\phi(z) = z + \frac{b_1}{z} + \frac{b_2}{z^2} + \ldots$ and $P, Q$ two polynomials for which $\phi P = c\phi Q$. Then we also have $P = cQ$.*

*Proof.* Let $p, q$ be the degrees of $P$ and $Q$. Expanding the identity $\phi P = c\phi Q$ in power series around infinity we have

$$P(z) + O(\tfrac{1}{z})^p = cQ(z) + O(\tfrac{1}{z})^q,$$

so necessarily $p = q$ and $P(z) = cQ(z)$. □

The following Lemma can be found in [Bea1].

**Lemma 5.** *Let $P$ be a polynomial and $c \in \mathbb{C}$. Then $J_{cP} = J_P$ if and only if $c \in \Sigma$.*

**Proposition 6.** *If $f = az^n + \ldots$, $g = bz^n + \ldots$ have same degree $n \geq 2$ and $J = J_f = J_g$, then*
$$g = \tfrac{b}{a} f \qquad \text{and} \qquad \tfrac{b}{a} \in \Sigma.$$

*Proof.* Let $\phi$ be as in Prop. 1. We have $\phi f = a(\phi)^n$ and $\phi g = b(\phi)^n$. This implies $\phi g = \tfrac{b}{a} \phi f$, and by Lemmas 4 and 5 we then have $g = \tfrac{b}{a} f$, and $\tfrac{b}{a} \in \Sigma$. □

Theorem 1 in [Bea2] states that $\{ g \mid \deg(f) = \deg(g), J_f = J_g \} = \{ \sigma f \mid \sigma \in \Sigma \}$. In our notation this is $\mathcal{S}_n = \Sigma f$, which is a direct corollary of Lemma 5 and Prop. 6. This proof is substantially different from that in [Bea2]. Moreover, by considering $f^i$, which has degree $n^i$, one has

**Proposition 7.** $\mathcal{S}_{n^i} = \Sigma f^i \;\; \forall i \geq 1$.

Also as a corollary, if $f$ and $g$ have different degrees we consider the polynomials $gf$ and $fg$, which have leading coefficients $ba^m, ab^n$, and get

**Proposition 8.** *If $f = az^n + \ldots$, $g = bz^m + \ldots$ have degrees $n, m \geq 2$ and $J = J_f = J_g$, then*
$$fg = \tfrac{b^{n-1}}{a^{m-1}} gf \qquad \text{and} \qquad \tfrac{b^{n-1}}{a^{m-1}} \in \Sigma.$$

From theorem 5 in [Bea1], if $l < \infty$ is the order of the group of symmetries $\Sigma$ of a centered polynomial $f$, then $f$ can be written in the form $f(z) = z^r P(z^l)$ for some polynomial $P$ and positive integer $r$.

The following can also be found in [Bea1]. We skip the proof.

**Lemma 9.** *If $|\sigma| = 1$, $\sigma \in \Sigma \iff \phi(\sigma z) = \sigma \phi(z) \iff f(\sigma z) = \sigma^n f(z)$.*

§3 *Some classical results on special Julia sets and commuting polynomials.*

For polynomials, there are two smooth Julia sets: a circle and an interval. They are related to the families $\{z^n\}$ and $\{T_n\}$ as follows (see, for example, [Bea3]):

If $S^1$ is the Julia set of a polynomial $P$ of degree $d \geq 2$, then $P(z) = \alpha z^d$, where $|\alpha| = 1$. If $[-1, 1]$ is the Julia set of a polynomial $P$ of degree $d \geq 2$, then $P(z) = \pm T_d$, where $T_d$ is the Tchebycheff polynomial of degree $d$.

From the work of Julia in [J], if $P, Q$ are two commuting polynomials of degrees $n, m \geq 2$, then we have one of the following three cases:

(1) $\exists L$ linear map and $\sigma \in \Sigma$ such that $P = Lz^n L^{-1}$ and $Q = L\sigma z^m L^{-1}$; or
(2) $\exists L$ linear map and $\sigma \in \Sigma$ such that $P = L T_n L^{-1}$ and $Q = L\sigma T_m L^{-1}$, where $T_n, T_m$ are Tchebycheff polynomials; or
(3) $\exists$ positive integers $\nu, \mu$ such that $P^\nu = Q^\mu$ (common iterate.)

Ritt shows in [R1] that if $P, Q$ are two commuting polynomials, then we have either (1) or (2) above or

(3') $\exists L$ linear map, $\epsilon_1, \epsilon_2 \in \Sigma$, positive integers $i, j$ and a polynomial $G$ such that
$$P = L(\epsilon_1 G^i) L^{-1}, \qquad\qquad Q = L(\epsilon_2 G^j) L^{-1}.$$

In [R2] it is shown that if two polynomials have a common iterate, $P^\nu = Q^\mu$, then one also has (3').

## §4 Polynomials with the same Julia set. Proof of Theorem 1.

*Proof of Theorem 1.*

1) If $J$ is a circle, consider a linear map $L_1(z) = Az + B$ such that $L_1(S^1) = J$. We have then that the maps $L_1^{-1} f L_1$ and $L_1^{-1} g L_1$ have $S^1$ as Julia set. Therefore

$$L_1^{-1} f L_1 = \alpha z^n, \qquad L_1^{-1} g L_1 = \beta z^m,$$

for some constants $\alpha, \beta$ with $|\alpha| = |\beta| = 1$. By conjugating now by an appropriate linear map $z \mapsto \lambda z$ and taking $L(z) = L_1(\lambda z)$ we get

$$L^{-1} f L = z^n, \qquad L^{-1} g L = \sigma z^m,$$

for some constant $\sigma$ with $|\sigma| = 1$.

2) If $J$ is an interval, consider a linear map $L(z) = Az + B$ such that $L([-1, 1]) = J$. We have then that the maps $L^{-1} f L$ and $L^{-1} g L$ have the interval $[-1, 1]$ as Julia set. By changing the sign of $A$ if necessary, we have

$$L^{-1} f L = T_n, \qquad L^{-1} g L = \pm T_m.$$

3) Since $J$ is not a circle, $\Sigma$ is finite. Let $l$ be the order of $\Sigma$. As observed earlier, we can write

$$f = z^{r_1} f_0(z^l), \qquad g = z^{r_2} g_0(z^l),$$

for some polynomials $f_0, g_0$ and positive integers $r_1, r_2$. In [B-E] these polynomials are related to the following commuting ones (see also [R3]): define

$$\hat{f} = z^{r_1} [f_0(z)]^l, \qquad \hat{g} = z^{r_2} [g_0(z)]^l.$$

Let $\psi(z) = z^l$, we have

$$\psi f = \hat{f} \psi, \qquad \psi g = \hat{g} \psi.$$

In [Bea1] it is observed that the Julia sets satisfy

$$\psi(J) = J_{\hat{f}} = J_{\hat{g}}.$$

The polynomials $\hat{f}, \hat{g}$ commute. Indeed, by Prop. 7, $fg = \sigma gf$ with $\sigma \in \Sigma$. We therefore have that $\psi(\sigma) = \sigma^l = 1$ and

$$\hat{f} \hat{g} \psi = \hat{f} \psi g = \psi f g = \psi \sigma g f = \psi g f = \hat{g} \psi f = \hat{g} \hat{f} \psi.$$

Therefore $\hat{f} \hat{g} \psi = \hat{g} \hat{f} \psi$, and this implies $\hat{f} \hat{g} = \hat{g} \hat{f}$.

Now, since we are in the case in which the Julia set $J$ is neither a circle nor an interval, $\psi(J) = J_{\hat{f}} = J_{\hat{g}}$ is neither a circle nor an interval. Therefore the commuting polynomials $\hat{f}, \hat{g}$ have a common iterate and there are positive integers $\nu, \mu$ such that
$$\hat{f}^\nu = \hat{g}^\mu.$$
As $\deg(f) = \deg(\hat{f})$ and $\deg(g) = \deg(\hat{g})$, we have
$$n^\nu = m^\mu.$$
This implies that there are positive integers $r, s, S$ such that
$$n = r^s, \qquad m = r^S.$$
Let $\phi$ be the map guaranteed by Prop. 1. We have
$$f(z) = \phi^{-1}\{a[\phi(z)]^{r^s}\}, \qquad g(z) = \phi^{-1}\{b[\phi(z)]^{r^S}\}.$$

The following can be found in [R2], and for completeness we include the proof.

**Lemma** (Ritt). *If the functions*
$$U(z) = \phi^{-1}\{h[\phi(z)]^{r^u}\}, \qquad V(z) = \phi^{-1}\{k[\phi(z)]^{r^v}\},$$
*where $u \geq v$, are polynomials, then the function*
$$UV^{-1}(z) = \phi^{-1}\{hk^{-r^{u-v}}[\phi(z)]^{r^{u-v}}\},$$
*where any of the determinations of $V^{-1}(z)$ is taken, is also a polynomial.*

*Proof of Lemma.* From the expression of $UV^{-1}(z)$ in terms of $\phi(z)$, we see that $UV^{-1}(z)$ is uniform in a neighbourhood of infinity. Let its Laurent development be
$$G(z) + P(1/z),$$
where $G(z)$ is a polynomial, and $P(1/z)$ a series of negative powers of $z$. We have
$$U(z) = G[V(z)] + P[1/V(z)].$$

The Laurent development of $P[1/V(z)]$ at infinity contains only negative powers of $z$, and since $U(z)$ has no negative powers in its expansion, $P[1/V(z)]$, and hence also $P[1/z]$, must be identically zero. Therefore
$$UV^{-1}(z) = G(z),$$
which proves the lemma. $\square$

Again following Ritt, let $t$ be the smallest integer greater than zero such that there exists a polynomial of the form
$$R(z) = \phi^{-1}\{e[\phi(z)]^{r^t}\}.$$

Then $t$ must be a divisor of $s$. Otherwise we would have (by Euclidean division algorithm)
$$s = it + j,$$
where $i$ and $j$ are integers and $0 < j < t$, and by the previous lemma, the function $f[g^{-i}(z)]$ would be a polynomial of the form
$$\phi^{-1}\{e'[\phi(z)]^{r^j}\}.$$

Similarly $t$ is a divisor of $S$.

We now proceed as follows. Let $q = s/t, Q = S/t$. From the previous lemma,

(∗)  $$f \circ (R^q)^{-1}(z) = \phi^{-1}[\epsilon_1 \phi(z)]$$

is a polynomial of degree one, and $\epsilon_1$ a constant. We have then that
$$\epsilon_1 \phi(z) = \phi(Az + B),$$
and from the expansion of $\phi$, it is clear that $A = \epsilon_1$, $B = 0$. Therefore

(∗∗)  $$\epsilon_1 \phi(z) = \phi(\epsilon_1 z).$$

From Lemma 9, we then have that $\epsilon_1 \in \Sigma$, and from (∗) and (∗∗)
$$f(z) = \epsilon_1 R^q(z).$$

Similarly, for some $\epsilon_2 \in \Sigma$ we will have
$$g(z) = \epsilon_2 R^Q(z).$$

This completes the proof of Theorem 1.  □

§5 *Commuting pairs and polynomials with same Julia set. Proof of Theorem 2.*

As we have mentioned, for polynomials related to the examples (1) and (2) in the introduction, we have that $\mathcal{S}_m \neq \emptyset \neq \mathcal{C}_m \; \forall m \geq 2$. However, in general we have $\mathcal{S}_m = \mathcal{C}_m = \emptyset$ for most $m$. We now study this in detail.

*Proof of Theorem 2.*

It is clear that $f^i$ has degree $n^i$ and $J_{f^i} = J_f$, so $\mathcal{S}_{n^i} \neq \emptyset$. Conversely, suppose $\mathcal{S}_m \neq \emptyset$ and let $g \in \mathcal{S}_m$. By Theorem 1, and as $f$ is minimal, there exists a polynomial $R$ and a positive integer $i$ such that
$$f = \epsilon_1 R, \qquad g = \epsilon_2 R^i.$$
We have $m = \deg(g) = (\deg R)^i = (\deg f)^i = n^i$. This proves (1).

As we already mentioned, we get (2) as a direct corollary of Propositions 5 and 6.

To prove (3), note that $\mathcal{S} = \cup \mathcal{S}_m$, and by (2), $\mathcal{S}_{n^i} = \Sigma f^i \subset \Sigma \mathcal{C}$. So $\mathcal{S} \subset \Sigma \mathcal{C}$. Conversely, let $g \in \mathcal{C}$ and $\sigma \in \Sigma$. As $g$ commutes with $f$, $J_g = J_f$ and, by Lemma 5, $J_{\sigma g} = J_g$. Therefore $\sigma g \in \mathcal{S}$.  □

PAU ATELA, DEPARTMENT OF MATHEMATICS, SMITH COLLEGE (NORTHAMPTON, MA 01063) & MSRI
*E-mail address*: patela@math.smith.edu

JUN HU, DEPARTMENT OF MATHEMATICS, CUNY & MSRI
*E-mail address*: huj@cunyvms1.gc.cuny.edu